\documentclass[11pt]{article} % use larger type; default would be 10pt

\usepackage[utf8]{inputenc} % set input encoding (not needed with XeLaTeX)

\usepackage{geometry} % to change the page dimensions
\usepackage{graphicx} % support the \includegraphics command and options

\usepackage{booktabs} % for much better looking tables
\usepackage{array} % for better arrays (eg matrices) in maths
\usepackage{paralist} % very flexible & customisable lists (eg. enumerate/itemize, etc.)
\usepackage{verbatim} % adds environment for commenting out blocks of text & for better verbatim
\usepackage{subfig} % make it possible to include more than one captioned figure/table in a single float
\usepackage[all]{xy}
\usepackage{amsmath} 
\usepackage{amsthm}
\usepackage{amsfonts}
\usepackage{fancyhdr} % This should be set AFTER setting up the page geometry
\usepackage{sectsty}
\allsectionsfont{\sffamily\mdseries\upshape} % (See the fntguide.pdf for font help)
\usepackage[nottoc,notlof,notlot]{tocbibind} % Put the bibliography in the ToC
\usepackage[titles,subfigure]{tocloft} % Alter the style of the Table of Contents

\title{On the existence of a torsor structure for Galois covers over a complete discrete valuation ring}
\author{Mohamed Sa\"\i di \& Nicholas Williams}
\date{} % Activate to display a given date or no date (if empty),
\newcommand{\Fr}{\mathrm{Fr}}
\newcommand{\defeq}{:=}
\newcommand{\Spec}{\mathrm{Spec}}
\newcommand{\chr}{\mathrm{char}}

\newtheorem{question}{Question}
\newtheorem*{thmA}{Theorem A}
\newtheorem*{thmB}{Theorem B}
\newtheorem*{thmC}{Theorem C}
\newtheorem*{remarksD}{Remarks D}

\begin{document}
\maketitle

\begin{abstract} In this note we investigate the problem of existence of a torsor structure for Galois covers 
of (formal) schemes over a complete discrete valuation ring of residue characteristic $p>0$ in the case of   
abelian Galois $p$-groups.
\end{abstract}

\section*{\S 0. Introduction}

In this paper $R$ denotes a \emph{complete discrete valuation ring}, with uniformiser $\pi$, residue field $k$ of characteristic $p>0$, 
and fraction field $K := \Fr R$.
For an $R$-(formal)scheme $Z$ we write $Z_K\defeq Z\times _{\Spec R}\Spec K$ and $Z_k\defeq Z\times _{\Spec R}\Spec k$ for 
the generic and special fibre, respectively, of $Z$. (In the case where $Z$ is a formal $R$-scheme 
by its generic fibre $Z_K$ we mean the associated rigid analytic space.) 
Let $X$ be a (\emph{formal}) $R$-scheme of finite type which is \emph{normal}, geometrically connected, and flat over $R$. We further assume 
that the special fibre $X_k$ of $X$ is \emph{integral}.
Let $f_K:Y_K\rightarrow X_K$ be an \emph{\'etale torsor} under a finite \'etale $K$-group scheme $\widetilde G$ of rank $p^t$ ($t\ge 1$), 
with $Y_K$ \emph{geometrically connected}, and $f:Y\rightarrow X$ the corresponding 
morphism of \emph{normalisation}. (Thus, $Y$ is the normalisation of $X$ in $Y_K$.)
We are interested in the following question.

\begin{question} When is $f:Y\rightarrow X$ a \emph{torsor} under a finite and flat $R$-group scheme $G$ which extends $\widetilde G$, i.e., with $G_K=\widetilde G$? 
\end{question}

The following is well known.

\begin{thmA} \emph{(Proposition 2.4 in [Sa\"\i di]; Theorem 5.1 in [Tossici])} If $\chr(K)=0$ we assume that $R$ contains a primitive $p$-th root of $1$, and $X$ is \emph{locally factorial}. Let $\eta$ be the generic point of $X_k$ and ${\mathcal{O}}_{\eta}$ the local ring of $X$ at $\eta$, which is a discrete valuation ring with fraction field $K(X)$: the function field of $X$. Let $f_K : Y_K \rightarrow X_K$ be an \'etale torsor under a finite \'etale $K$-group scheme $\widetilde G$ of {\bf rank} $\bold p$, with $Y_K$ connected, and let $K(X) \rightarrow L$ be the corresponding extension of function fields. Assume that the ramification index above ${\mathcal{O}}_\eta$ in the field extension  $K(X) \rightarrow L$ equals 1. Then $f : Y \rightarrow X$ is a torsor under a finite and flat $R$-group scheme $G$ of rank $p$ which extends $\widetilde{G}$ (i.e., with $G_K=\widetilde{G}$). 
\end{thmA}

Strictly speaking the above references treat the case where $\chr(K)=0$. For the equal characteristic $p>0$ case see
[Sa\"\i di1], Theorem 2.2.1. Theorem A also holds when $X$ is the formal spectrum of a complete discrete valuation ring (cf. 
[Sa\"\i di2], Proposition 2.3, and the references therein in the unequal characteristic case, as well as Proposition 2.3.1 in 
[Sa\"\i di3] in the equal characteristic $p>0$ case). 
It is well known that the analog of Theorem A is \emph{false} in general. 
There are counterexamples to the statement in Theorem A where $\widetilde G$ is cyclic of rank $p^2$, 
see [Tossici], Example 6.2.12, for instance. 

Next, we describe the setting in this paper. Let $n\ge 1$, and for  $i\in \{1,\cdots,n\}$ let 
$$f_{i,K}:X_{i,K} \rightarrow X_K$$ 
be an \'etale torsor under an \'etale finite \emph{commutative} $K$-group scheme $\widetilde G_i$, with $X_{i,K}$ \emph{geometrically connected},
such that the $\{f_{i,K}\}_{i=1}^n$ are \emph{generically pairwise disjoint}. 
Assume that $f_{i,K}:X_{i,K} \rightarrow X_K$ extends to a torsor
$$f_i:X_i\rightarrow X$$ 
under a finite and flat (necessarily commutative) $R$-group scheme $G_i$ with $(G_i)_K=\widetilde G_i$, 
and with $X_i$ \emph{normal}, $\forall i\in \{1,\cdots,n\}$. (Thus, $X_i$ is the normalisation of $X$ in $X_{i,K}$.)
Let 
$$\widetilde {X}_K\defeq X_{1,K}\times _{X_K}X_{2,K}\times _{X_K} \cdots \times _{X_K}X_{n,K},$$ 
and $\widetilde X$ the \emph{normalisation} of $X$ in $\widetilde X_K$. Thus, $\widetilde X_K$ is the generic fibre
of $\widetilde X$ and we have the following commutative diagrams

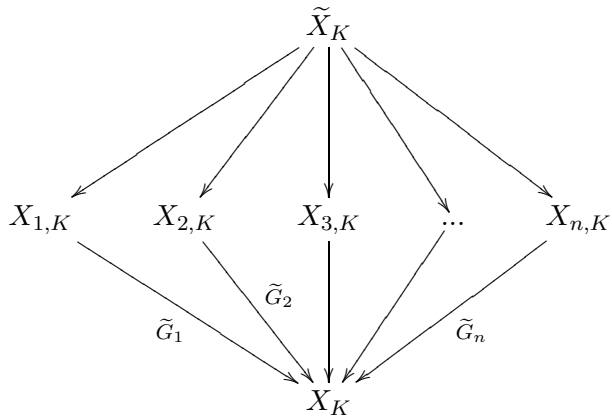
\begin{equation*}
 \xymatrix{
 && \widetilde{X}_K \ar@{->}[ddll] \ar@{->}[ddl] \ar@{->}[dd] \ar@{->}[ddr] \ar@{->}[ddrr] \\
  \\
X_{1,K} \ar@{->}[ddrr]_{\widetilde{G}_{1}} & X_{2,K} \ar@{->}[ddr]^{\widetilde{G}_{2}}  & X_{3,K} \ar@{->}[dd] & ... \ar@{->}[ddl] & X_{n,K} \ar@{->}[ddll]^{\widetilde{G}_{n}} \\
  \\
&& X_K 
}
\end{equation*}

and
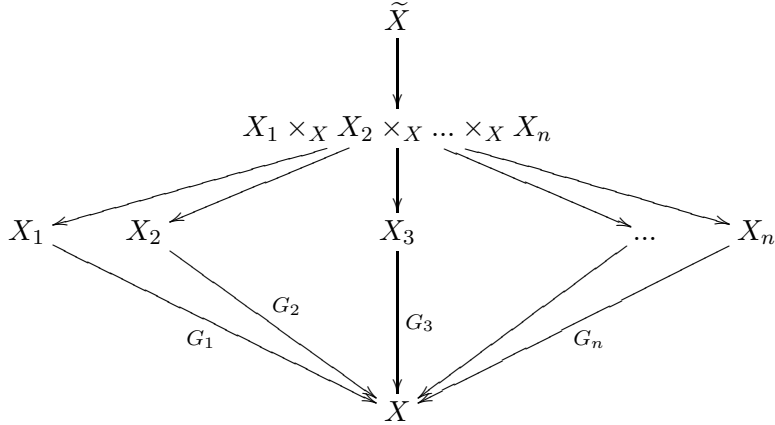
\begin{equation*}
 \xymatrix{
&& \widetilde{X} \ar@{->}[d]  \\
&& X_{1} \times_X X_{2} \times_X  ... \times_X X_{n} \ar@{->}[dll] \ar@{->}[dl] \ar@{->}[d] \ar@{->}[dr] \ar@{->}[drr]  \\
X_{1} \ar@{->}[ddrr]_{G_{1}} & X_{2} \ar@{->}[ddr]^{G_{2}}  & X_{3} \ar@{->}[dd]^{G_{3}} & ...  \ar@{->}[ddl] & X_{n} \ar@{->}[ddll]^{G_{n}} \\
\\
&& X
}
\end{equation*}
where $X_{1} \times_X X_{2} \times_X \cdots \times_X X_{n}$ denotes the fibre product of the $\{X_i\}_{i=1}^n$ over $X$, the morphism 
$\widetilde X\rightarrow X_{1} \times_X X_{2} \times_X \cdots \times_X X_{n}$ is birational and is induced by the natural \emph{ finite} morphisms
$\widetilde X\rightarrow X_{i}$, $\forall i\in \{1,\cdots,n\}$. Note that 
$f_K:\widetilde X_K\rightarrow X_K$ (resp. $\tilde f:X_{1} \times_X X_{2} \times_X \cdots \times_X X_{n} \rightarrow X$) is a torsor 
under the \'etale finite commutative $K$-group scheme $\widetilde G\defeq \widetilde G_1 \times_{\Spec K} \widetilde G_2 \times_{\Spec K} \cdots 
\times_{\Spec K} \widetilde G_n$
(resp. a torsor under the finite and flat commutative $R$-group scheme $G_1 \times_{\Spec R} G_2 \times_{\Spec R} \cdots \times_{\Spec R} G_n$), 
as follows easily from the various definitions. Note that $\left (G_1 \times_{\Spec R} G_2 \times_{\Spec R} \cdots \times_{\Spec R} G_n 
\right )_K=\widetilde G$. 

In this setup Question 1 reads as follows.

\begin{question} When is $f:\widetilde{X} \rightarrow X$ a torsor under a finite and flat (necessarily commutative) $R$-group scheme 
$G$ which extends $\widetilde G$, i.e., with $G_K=\widetilde G$?
\end{question}

Our main result in this paper is the following. 

\begin{thmB} We use the same notations as above. Assume that $\widetilde{X}_k$ is {\bf reduced}. Then the following three statements are equivalent.
\begin{enumerate}
\item $f:\widetilde{X}\rightarrow X$ is a torsor under a finite and flat commutative $R$-group scheme $G$, in which case  $G=G_1 \times_{\Spec R} \cdots \times_{\Spec R} G_n$ necessarily.
\item $\widetilde{X}= X_{1} \times_X X_{2} \times_X  \cdots \times_X X_{n}$, in other words 
$X_{1} \times_X X_{2} \times_X  \cdots \times_X X_{n}$ is normal.
\item $\left( X_{1} \times_X X_{2} \times_X \cdots\times_X X_{n} \right)_k$ is reduced.
\end{enumerate}
\end{thmB}

Note that the above condition in Theorem B that $\widetilde{X}_k$ is reduced is always satisfied after possibly passing to a finite extension $R'/R$ of $R$ (cf. [Epp]).  It implies that the $(X_i)_k$ are reduced, $\forall i\in \{1,\cdots,n\}$. Moreover, Theorem A and Theorem B provide a ``complete" answer to Question 1 in the case of Galois covers of type $\left (p,\cdots,p\right )$, i.e., the case where $\mathrm{rank} (G_i)=p,\ \forall i\in \{1,\cdots,n\}$.
 
In the case of (relative) \emph{smooth curves} one can prove the following more precise result.

\begin{thmC}  We use the same notations and assumptions as in Theorem B. Assume further that $X$ is a (relative) \textbf{smooth $R$-curve}, 
$n\ge 2$, and $R$ is \textbf{strictly henselian}. If $\chr(K)=0$ we assume that $K$ contains a primitive $p$-th root of $1$.
Then the three (equivalent) conditions in Theorem B are equivalent to the following.

4. \textbf{At least} $\textbf{n-1}$ of the finite flat $R$-group schemes $G_{i}$ acting on $f_i:X_i \rightarrow X$ are \textbf{ \'{e}tale}, for $i \in \{1,\cdots,n\}$.
\end{thmC}

\begin{remarksD}
\begin{enumerate}
\item Theorem B holds true if $X$ is the formal spectrum of a complete discrete valuation ring 
(cf. the details of the proof of Theorem B in $\S1$ which applies as it is in this case).
\item In $\S3$ we provide examples showing that Theorem C \emph{doesn't hold} in relative dimension $>1$.
\end{enumerate}
\end{remarksD}

\section*{\S 1. Proof of Theorem B}

In this section we prove Theorem B. We start by the following.

\textbf{Proposition 1.1}
\emph{Let $G$ be a finite and flat commutative $R$-group scheme whose generic fibre is a product of group schemes of 
the form
$$ G_K=\widetilde G_1 \times_{\Spec K}\widetilde G_2 \cdots  \times_{\Spec K} \widetilde G_n,$$
where the $\{\widetilde G_i\}_{i=1}^n$ are finite and flat commutative $K$-group schemes. Then $G$ is a product of 
finite and flat commutative $R$-group schemes  $\{G_i\}_{i=1}^n$, i.e.,
$$ G = G_1 \times_{\Spec R}  G_2 \times_{\Spec R}\cdots \times_{\Spec R} G_n,$$
with $(G_{i})_K= \widetilde G_i$.}

\begin{proof} First, we treat the case $n=2$. Thus, we have $G_K=\widetilde G_1 \times_{\Spec K} \widetilde G_2$ and need 
to show $G = G_1 \times_{\Spec R} G_2$ where $(G_{i})_K= \widetilde G_i$, for $i=1,2$. 
Let $G_i$ be the \emph{ schematic closure} of $\widetilde G_{i}$ in $G$, for $i=1,2$ (cf. [Raynaud], 2.1). Therefore, $G_1$ and $G_2$ are 
closed subgroup schemes of $G$ which are finite and flat over $\Spec R$ (cf. loc. cit.). 
We have a short exact sequence
$$ 1  \rightarrow  G_1  \rightarrow  G  \rightarrow  G/G_1  \rightarrow 1 ,$$
and likewise
$$ 1  \rightarrow  G_2  \rightarrow  G  \rightarrow  G/G_2  \rightarrow 1 ,$$
of finite and flat commutative $R$-group schemes (cf. loc. cit.).
It remains for the proof to show that the composite homomorphism $G_2 \rightarrow G \rightarrow G/G_1$ is 
an isomorphism. The morphism $G \rightarrow G/G_1$ is finite. The morphism $G_2 \rightarrow G$ is a closed 
immersion, hence finite. The composite $G_2 \rightarrow G/G_1$ of the above morphisms is then finite. 
We will show it is an isomorphism. The morphism $G_2 \rightarrow G/G_1$ 
is a closed immersion since its kernel is trivial. Indeed, on the generic fibre the kernel is trivial: 
$\left(G_1 \cap G_2\right)_K = \widetilde{G_1} \cap \widetilde{G_2}  = \{1\}$. The map $G_2 \rightarrow G/G_1$ is then an 
isomorphism as both group schemes have the same rank. Similarly, the morphism $G_1 \rightarrow G/G_2$ is an isomorphism. 
Therefore, $G=G_1 \times_{\Spec R} G_2$ as required. Now an easy devissage argument along the above lines of thought, using 
induction on $n$, reduces immediately to the above case $n=2$.
\end{proof}

\emph{Proof of Theorem B}

\begin{proof}
(1 $\Rightarrow$ 2)
Assume that $f:\widetilde{X}\rightarrow X$ is a torsor under a finite and flat $R$-group scheme $G$.
In particular, $G_K=\widetilde G$ and $G$ is necessarily commutative.
We will show that $\widetilde{X}= X_{1} \times_X X_{2} \times_X  ... \times_X X_{n}$, i.e., show that 
$X_{1} \times_X X_{2} \times_X  ... \times_X X_{n}$ is normal (this will imply that  $G=G_1 \times_{\Spec R}\cdots \times_{\Spec R} G_n$
necessarily, as $G_1 \times_{\Spec R} ... \times_{\Spec R} G_n$ is the group scheme of the torsor $\tilde f:X_{1} \times_X X_{2} \times_X\cdots 
\times_X X_{n}\rightarrow X$). One reduces easily by a devissage argument to the case $n=2$ 
which we will treat below.

Assume $n=2$. We have the following commutative diagrams of torsors

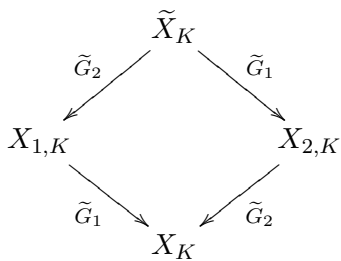
\begin{equation*}
 \xymatrix{
& \widetilde{X}_K  \ar@{->}[dl]_{\widetilde{G}_{2}} \ar@{->}[dr]^{\widetilde{G}_{1}}  \\
X_{1,K} \ar@{->}[dr]_{\widetilde{G}_{1}}  &&  X_{2,K}  \ar@{->}[dl]^{\widetilde{G}_{2}}   \\
& X_K 
}
\end{equation*}

and

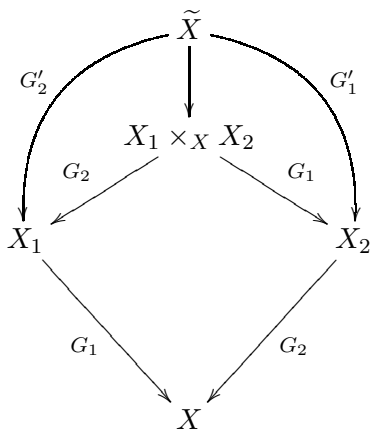
\begin{equation*}
 \xymatrix{
& \widetilde{X} \ar@{->}[d]   \ar@/^2pc/[ddr]^{G'_1}  \ar@/_2pc/[ddl]_{G'_2}   \\
& X_{1} \times_X X_{2}  \ar@{->}[dl]_{G_2  \ \ } \ar@{->}[dr]^{ \ \ G_1}  \\
X_{1} \ar@{->}[ddr]_{G_1}  &&  X_2  \ar@{->}[ddl]^{G_2}   \\
\\
& X 
}
\end{equation*}
where $\widetilde X\rightarrow X_i$ is a torsor under a finite and flat $R$-group scheme $G_i'$, for $i=1,2$.
Moreover, $G'_1 = \left( \widetilde G_{1} \right)^{\text{schematic closure}}$, and $G'_2 = \left( \widetilde G_{2} 
\right)^{\text{schematic closure}}$ (where the schematic closure is taken inside $G$) holds necessarily, so that $G=G'_1 \times_{\Spec R} G'_2$ (cf. Proposition 1.1).
Note that $\widetilde X/{G'_1}=X_2$ must hold as the quotient $\widetilde X/{G'_1}$ is normal: since  
$\left (\widetilde X/{G'_1}\right )_k$ is reduced (as $\widetilde X_k$ is reduced and $\widetilde X$ dominates $\widetilde X/{G'_1}$),
and $\left (\widetilde X/{G'_1}\right )_K=X_{2,K}$
is normal (cf. [Liu], 4.1.18). Similarly $\widetilde X/{G'_2}=X_1$ holds.
We want to show that $\widetilde{X}= X_{1} \times_X X_{2}$, and we claim that this reduces to showing that the natural morphism $G \to G_1 \times_{\Spec R} G_2$ 
(cf. the map $\phi$ below) is an isomorphism. Indeed, if one has two torsors, in this case $ \widetilde{X}\rightarrow X$
and $X_{1} \times_X X_{2}\rightarrow X$ above the same $X$, under isomorphic group schemes, which are isomorphic on the generic fibres, and if we have a morphism 
$\widetilde{X} \rightarrow X_{1} \times_X X_{2}$ which is compatible with the torsor structure and the given identification of group schemes (cf. above diagrams and the definition of $\phi$ below), 
then this morphism must be an isomorphism. (This is a consequence of Lemma 4.1.2 in
[Tossici]. In [Tossici] $\chr(K)=0$ is assumed, the same proof however applies if $\chr(K)=p$.) 
We have two short exact sequences of finite and flat commutative $R$-group schemes
(cf. above diagrams and discussion for the equalities  $G_1=G/G'_2$ and  $G_2=G/G'_1$)
$$ 1  \rightarrow G'_2  \rightarrow G  \rightarrow G_1 = G/G'_2  \rightarrow 1 ,$$
and 
$$ 1  \rightarrow G'_1  \rightarrow G  \rightarrow G_2 = G/G'_1  \rightarrow 1 .$$

The morphisms $G \rightarrow G_1$, and $G \rightarrow G_2$, are finite. Consider the following exact sequence

$$ 1  \rightarrow \mathrm{Ker}(\phi) \rightarrow G \rightarrow G_1 \times_{\Spec R} G_2 ,$$
where $\phi:G   \rightarrow G_1 \times_{\Spec R} G_2$ is the morphism induced by the above morphisms.
We want to show that the map $\phi:G \rightarrow G_1 \times_{\Spec R} G_2$ is an isomorphism. We have $\mathrm{Ker}(\phi)=G'_1 \cap G'_2$ 
by construction. However, $G'_1 \cap G'_2 = \{ 1 \}$ since $G=G'_1 \times _{\Spec R}G'_2$ by Proposition 1.1, and therefore 
$\mathrm{Ker}(\phi) = \{ 1 \}$ which means $\phi:G \rightarrow G_1 \times_{\Spec R} G_2$  is a closed immersion. Finally, $G$ and 
$G_1 \times_{\Spec R} G_2$ have the same rank as group schemes which implies $\phi$ is an isomorphism, as required. 

(2 $\Rightarrow$ 3) Clear.

(3 $\Rightarrow$ 1) By assumption $\left(X_{1} \times_X X_{2} \times_X  ... \times_X X_{n} \right)_k$ is reduced. 
Moreover, we have $\left( X_{1} \times_X X_{2} \times_X  ... \times_X X_{n} \right)_K=\widetilde X_K$ is normal. Hence 
$X_{1} \times_X X_{2} \times_X  ... \times_X X_{n}$ is normal (cf. [Liu], 4.1.18),
and $\widetilde{X} = X_{1} \times_X X_{2} \times_X  ... \times_X X_{n}$. We know 
that $\tilde f:X_{1} \times_X X_{2} \times_X  ... \times_X X_{n}\rightarrow X$ 
is a torsor under the group scheme 
$G_1 \times_{\Spec R} G_2 \times_{\Spec R} .... \times_{\Spec R} G_n$, 
so $f:\widetilde{X}\rightarrow X$ is a torsor under the group scheme $G=G_1 \times_{\Spec R}  G_2 \times_{\Spec R}  .... \times_{\Spec R} G_n$.
\end{proof}

\section*{\S 2. Proof of Theorem C }
In this section we prove Theorem C.

\begin{proof}

(1 $\Rightarrow$ 4) Suppose that $\tilde f:\widetilde{X}\rightarrow X$ is a torsor under a finite and flat $R$-group scheme $G$; 
in which case $\widetilde{X} = X_{1} \times_X X_{2}\times_X  ... \times_X X_{n}$ and $G=G_1 \times_{\Spec R} \cdots \times_{\Spec R} G_n$ 
(cf. Theorem B). We will show that \emph{at least} ${n-1}$ of the finite 
flat $R$-group schemes $G_{i}$ (acting on $f_i:X_i \rightarrow X$) are \'{e}tale, for $i\in \{1,\cdots,n\}$.
We argue by induction on the rank of $G$.

\emph{ Base case:} The base case pertains to  $\mathrm{rank} (G)=p^2$ and $n=2$. Thus, $\mathrm{rank} (G_1)=\mathrm{rank} (G_2)=p$.
We assume $\widetilde{X} = X_{1} \times_X X_{2}$ and prove that at 
least one of the two group schemes $G_1$ or $G_2$ is \'{e}tale. We assume that $X$ is a scheme, and not a formal scheme, in which case the argument of proof is the same.

Let $x$ be a \emph{closed} point of $X$ and $\mathcal{X}$ the \emph{boundary of the formal germ} of $X$ at $x$, so $\mathcal{X}$  is isomorphic to 
$\mathrm{Spec} \left( R[[T]] \{ T^{-1} \} \right)$ (cf. [Sa\"\i di2], \S 1). We have a natural morphism $\mathcal{X}\rightarrow X$ of schemes.
Write $\mathcal{X}_1\defeq \mathcal{X}\times _XX_1$, $\mathcal{X}_2\defeq \mathcal{X}\times _XX_2$, and $\widetilde {\mathcal{X}}\defeq 
\mathcal{X}\times _X\widetilde X$. Thus, by base change, $\widetilde {\mathcal{X}}\rightarrow \mathcal{X}$ (resp. $\mathcal{X}_1\rightarrow \mathcal{X}$, and 
$\mathcal{X}_2\rightarrow \mathcal{X}$) is a torsor under the group scheme $G$ (resp. under $G_1$, and $G_2$) and we have the following commutative diagram 

\begin{equation*}
 \xymatrix{
& \widetilde{\mathcal{X}} = \mathcal{X}_{1} \times_\mathcal{X} \mathcal{X}_{2}  \ar@{->}[dl]_{G_{2}} \ar@{->}[dr]^{G_{1}}  \\
\mathcal{X}_{1} \ar@{->}[dr]_{G_{1}}  &&  \mathcal{X}_2  \ar@{->}[dl]^{G_{2}}   \\
& \mathcal{X}
}
\end{equation*}
Note that 
$\widetilde {\mathcal{X}}$ is normal as $(\widetilde {\mathcal{X}})_k$ is reduced (recall $(\widetilde X)_k$ 
is reduced) and $(\widetilde {\mathcal{X}})_K$ 
is normal (cf. [Liu], 4.1.18), hence $\widetilde {\mathcal{X}}=\mathcal{X}_2\times _{\mathcal{X}}\mathcal{X}_2$ holds 
(cf. Theorem B and Remarks D, 1).

Assume now that $G_{1}$ and $G_{2}$ are both \emph{non-\'{e}tale} $R$-group schemes. Then we prove that  
$\widetilde {\mathcal{X}}\rightarrow \mathcal{X}$ can not have the structure of a torsor under a finite and flat $R$-group scheme which 
would then be a contradiction. More precisely, we will prove that $\mathcal{X}_2\times _{\mathcal{X}}\mathcal{X}_2$ can not be normal in this case, 
hence the above conclusion (cf. Theorem B).

We will assume for simplicity that $\chr(K)=0$ and $K$ contains a primitive $p$-th root of $1$.
A similar argument used below holds in equal characteristic $p>0$.
First, $\widetilde {\mathcal{X}}$ is connected as $\widetilde X_k$ is \emph{unibranch} (the finite morphism $\widetilde X_k\to X_k$ is radicial).
As the group schemes $G_1$ and $G_2$ are non  \'{e}tale, their special fibres $(G_{1})_k$ and $(G_{2})_k$ are radicial isomorphic to either 
$\mu_p$ or $\alpha_p$. We treat the case  $(G_{1})_k$ is isomorphic to $\mu_p\defeq \mu_{p,k}$ and $(G_{2})_k$ is isomorphic to $\alpha_p\defeq \alpha_{p,k}$, 
the remaining cases are treated similarly. 
(Recall $\mathcal{X}$ is isomorphic to $\mathrm{Spec} \left(R[[T]] \{ T^{-1} \} \right)$.) For a suitable choice of the parameter $T$ 
the torsor $\mathcal{X}_2\rightarrow \mathcal{X}$ is given by an equation $Z_2^p=1+\pi^{np}T^m$
where $n$ is a positive integer (satisfying a certain condition) and $m\in \mathbb{Z}$ (cf. [Sa\"\i di2], Proposition 2.3 (b). Strictly speaking in loc. cit.
this is shown to hold after a finite extension of $R$, however a close inspection of the proof in loc. cit. reveals that this finite extension can be chosen to be \'etale.
Also see Proposition 2.3.1 in [Sa\"\i di3] for the equal characteristic case),
and the torsor $\mathcal{X}_1\rightarrow \mathcal{X}$ is given by an equation $Z_1^p=f(T)$ where $f(T)\in R[[T]] \{ T^{-1} \}$ is a unit whose 
reduction $\overline {f(T)}$ modulo $\pi$ is not a $p$-power (cf. loc. cit.). We claim that 
$\widetilde{\mathcal{X}} = \mathcal{X}_{1} \times_{\mathcal{X}} {\mathcal{X}}_{2}$ can not hold. Indeed, by base change ${\mathcal{X}}_{1} \times_{\mathcal{X}} \mathcal{X}_{2}\rightarrow 
\mathcal{X}_2$ is a $G_1$-torsor 
which is generically given by an equation  $Z^p=f(T)$, where $f(T)$ is viewed as a function on $\mathcal{X}_2$. 
But in $\mathcal{X}_2$ the function $T$ becomes a $p$-power modulo $\pi$
as one easily deduces from the equation $Z_2^p=1+\pi^{np}T^m$ defining the torsor $\mathcal{X}_2\rightarrow \mathcal{X}$. In particular, the   
reduction $\overline {f(T)}$ modulo $\pi$ of $f(T)$, viewed as a function on $\left (\mathcal{X}_2\right )_k$, is a $p$-power. This means that 
$({\mathcal{X}}_{1} \times_{\mathcal{X}} \mathcal{X}_{2})_k$ is not reduced and
$\widetilde {\mathcal{X}}\rightarrow \mathcal{X}_2$ can not be a $G_1\simeq \mu_{p,R}$-torsor (cf. the proof of Proposition 2.3 in [Sa\"\i di2]), and a fortiori
$\widetilde{\mathcal{X}} \neq {\mathcal{X}}_{1} \times_{\mathcal{X}} \mathcal{X}_{2}$.

\emph{ Inductive hypothesis: } Given $G$, we assume that the (1 $\Rightarrow$ 4) part in Theorem C holds true for $n \ge 2$ and cases where
$ \mathrm{rank}(G_1)+\cdots+ \mathrm{rank} (G_n)< \mathrm{rank} (G)$. Write $\widetilde{X}_1 \defeq X_{1} \times_X X_{2} \times_X  ... \times_X X_{n-1}$. Then $\widetilde{X}_1$ 
is normal (since its special fibre is reduced (as it is dominated by $\widetilde{X}$ whose special fibre is reduced)
and its generic fibre is normal (cf. [Liu], 4.1.18)), hence at least $n-2$ of the 
corresponding $G_{i}$'s, for $i\in \{1,\cdots,n-1\}$, are \'{e}tale by the induction hypothesis. We will assume, without loss of generality, 
that $G_{i}$ is \'etale for $1 \leq i \leq n-2$.

\emph{ Inductive step: } We have the following picture for our inductive step (the case for $n$):

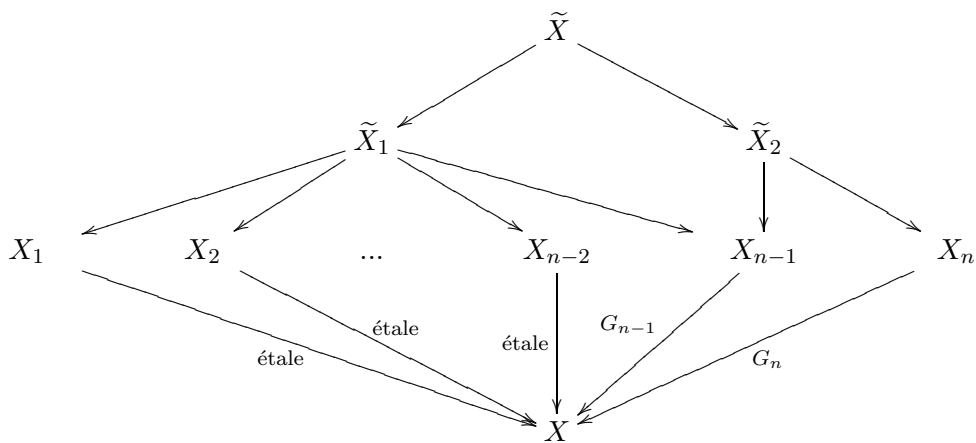
\begin{equation*}
 \xymatrix{
&&& \widetilde{X}\\
&&  \widetilde{X}_1  \ar@{<-}[ur]  && \widetilde{X}_2  \ar@{<-}[ul]  \\
\ \ \ X_{1} \ \ \ \ar@{<-}[urr]  & \ \ \  X_{2} \ \ \  \ar@{<-}[ur]  & \ \ \ ... \ \ \ & \ \ \ X_{n-2} \ \ \ \ar@{<-}[ul]  & \ \ \ X_{n-1} \ \ \ \ar@{<-}[ull] \ar@{<-}[u] & \ \ \ X_n \ \ \  \ar@{<-}[ul] \\
\\
&&& X   \ar@{<-}[uulll]^{\text{\'etale}}  \ar@{<-}[uull]_{\text{\'etale}} \ar@{<-}[uu]^{\text{\'etale}}  \ar@{<-}[uur]^{G_{n-1}}   \ar@{<-}[uurr]_{G_{n}}
}
\end{equation*}

We argue by contradiction. Suppose that neither $G_{n-1}$ nor $G_{n}$ is \'{e}tale.  This would mean that $\widetilde{X}_2\rightarrow X$, where $\widetilde{X}_2$ is
the normalisation of $X$ in $(X_{n-1})_K \times_{X_K} (X_{n})_K$, does not have the structure of a torsor (as this would contradict the induction hypothesis).
This implies that $\widetilde{X}\rightarrow X$ does not have the structure of a torsor since it factorises 
$\widetilde{X}\rightarrow \widetilde{X}_2\rightarrow X$, for otherwise $\widetilde{X}_2\rightarrow X$ being a quotient of 
$\widetilde{X}\rightarrow X$ would be a torsor. 
Of course, $\widetilde{X}\rightarrow X$ is a torsor to start with by assumption and so this is a contradiction. 
Therefore, at least one of $G_{n-1}$ and $G_{n}$ is \'{e}tale, as required.

(1 $\Leftarrow$ 4)  Suppose that at least $n-1$ of the $G_{i}$ are \'{e}tale, say: $G_1,G_2,\cdots,G_{n-1}$ are \'etale. Write 
$\widetilde X_1\defeq  X_{1} \times_X X_{2}\times_X  ... \times_X X_{n-1}$. Then $\widetilde X_1\rightarrow X$ is a torsor under the 
finite {\it \'etale} $R$-group scheme 
$G_1'\defeq G_1\times _{\Spec R}G_2\times_{\Spec R}\cdots\times _{\Spec R}G_{n-1}$. Moreover, $X_{1} \times_X X_{2}\times_X  ... \times_X X_{n}
=\widetilde X_1\times_XX_n$, 
and $X_{1} \times_X X_{2}\times_X  ... \times_X X_{n}\rightarrow X_{n}$ is an \'etale torsor under the group scheme $G'_1$ (by base change). 
In particular, $\left (X_{1} \times_X X_{2}\times_X  ... 
\times_X X_{n}\right )_k$ is reduced as $(X_{n})_k$ is reduced.
(Indeed, $\widetilde X$ dominates $X_n$ and ${\widetilde X}_k$ is reduced.) 
Hence $\widetilde X=X_{1} \times_X X_{2}\times_X  ... \times_X X_{n}$ (cf. Theorem B) and 
$\widetilde X\rightarrow X$ is a torsor under the group scheme $G\defeq G_1\times _{\Spec R}G_2\times_{\Spec R}\cdots\times _{\Spec R}G_{n}$.
\end{proof}

\section*{ \S 3. Counterexample to Theorem C in higher dimensions}

Theorem C is not valid (under similar assumptions) for (formal) smooth $R$-schemes of relative dimension $\geq 2$. Here is a counterexample.
Assume $\chr(K)=0$ and $K$ contains a primitive $p$-th root of $1$. Let $X= \mathrm{Spf}(A)$ where $A\defeq R<T_1, T_2>$ is the free 
$R$-Tate algebra in the two variables $T_1$ and $T_2$. 
Let $G_1=G_2=\mu_p\defeq \mu_{p,R}$, neither being an \'etale $R$-group scheme. For $i=1,2$, consider the 
$G_i$-torsor $X_i\rightarrow X$ which is generically defined by the equation 
$$Z_i^p=T_i.$$ 
We have the following commutative diagram

\begin{equation*}
 \xymatrix{
&   X_1 \times X_2  \ar@{->}[ddl]^{(Z'_2)^p=T_2}_{\mu_p} \ar@{->}[ddr]_{(Z'_1)^p=T_1}^{\mu_p}  \\
\\
X_1   \ar@{->}[ddr]_{\mu_p}^{Z_1^p=T_1}  && X_2  \ar@{->}[ddl]^{\mu_p}_{Z_2^p=T_2} \\
\\
&  X= \mathrm{Spf} \left( R<T_1, T_2> \right) \\
}
\end{equation*}

The torsor $X_1\times_XX_2\rightarrow X_2$ is a $G_1=\mu_p$-torsor defined generically by the equation 
$$(Z'_1)^p=T_1$$ 
where $T_1$ is viewed as a function on $X_2$.
This function is not a $p$-power modulo $\pi$ as follows easily from the fact that the torsor $X_2\rightarrow X$ is defined 
generically by the equation $Z_2^p=T_2$. In particular, $X_1\times_XX_2\rightarrow X_2$ is a non trivial $\mu_p$-torsor, and
$(X_1\times_XX_2)_k\rightarrow (X_2)_k$ is a non trivial $\mu_{p,k}$-torsor. Hence
$(X_1\times_XX_2)_k$ is necessarily reduced (as $(X_2)_k$ is reduced since $(X_2)_k\to (X_1)_k$ is a non trivial $\mu_{p,k}$-torsor). 
Thus,  $X_1\times_XX_2$ is normal (cf. Theorem B) and 
$X_1\times_XX_2=\widetilde X$, where $\widetilde X$ is the normalisation of $X$ in 
$(X_1\times_XX_2)_K$, which contradicts the statement of Theorem C in this case.

$$\textbf{References.}$$

\noindent
[Epp] Epp, H.P., Eliminating wild ramification, Inventiones Mathematicae, Volume 19, 1973, pp. 235-249.

\noindent
[Liu] Liu, Q., Algebraic Geometry and Arithmetic Curves, Oxford University Press, 2009.

\noindent
[Raynaud] Raynaud, M., Sch\'emas en groupes de type $(p,...,p)$, Bulletin de la Soci\'et\'e Math\'ematique de France, 
Volume 102, 1974, pp. 241-280.

\noindent
[Sa\"\i di] Sa\"\i di, M., Torsors under finite and flat group schemes of rank $p$ with Galois action, Mathematische 
Zeitschrift, Volume 245, 2003, pp. 695-710.

\noindent
[Sa\"\i di1] Sa\"\i di, M., On the degeneration of \'etale $\Bbb Z/p\Bbb Z$ and $\Bbb Z/p^2\Bbb Z$-torsors in equal characteristic $p>0$,
Hiroshima  Math. J. 37 (2007), 315-341.

\noindent
[Sa\"\i di2] Sa\"\i di, M., Wild ramification and a vanishing cycles formula, Journal of Algebra 273 (2004) 108-128.

\noindent
[Sa\"\i di3] Sa\"\i di, M., Galois covers of degree $p$ and semi-stable reduction of curves in equal characteristic $p>0$,
Math. J. Okayama Univ. 49 (2007), 113-138.

\noindent
[Tossici] Tossici, D., Effective models and extension of torsors over a discrete valuation ring of unequal characteristic, International Mathematics Research Notices, 2008, pp. 1-68.

\bigskip

\noindent
Mohamed Sa\"\i di

\noindent
College of Engineering, Mathematics, and Physical Sciences

\noindent
University of Exeter

\noindent
Harrison Building

\noindent
North Park Road

\noindent
EXETER EX4 4QF

\noindent
United Kingdom

\noindent
M.Saidi@exeter.ac.uk

\bigskip

\noindent
Nicholas Williams

\noindent
London, United Kingdom

\noindent
Nicholas.Williams09@alumni.imperial.ac.uk

\end{document}